\documentclass[11pt]{amsart}
\usepackage{amsmath,amsthm,epsfig,amssymb}
\title{Apollonius ``circle" in Hyperbolic Geometry}
\author{Eugen J. Iona\c scu}
\curraddr{CSU, Math Dept, Columbus, GA, USA } \email{math@ejionascu.ro} \subjclass{51N20}
\date{December $5^{nd}$, 2016}
    \keywords{Hyperbolic Space, Circle of Apollonius, polar coordinates, geometric probability }

\begin{document}

\newtheorem{theorem}{\hspace{\parindent}
T{\scriptsize HEOREM}}[section]
\newtheorem{proposition}[theorem]
{\hspace{\parindent }P{\scriptsize ROPOSITION}}
\newtheorem{corollary}[theorem]
{\hspace{\parindent }C{\scriptsize OROLLARY}}
\newtheorem{lemma}[theorem]
{\hspace{\parindent }L{\scriptsize EMMA}}
\newtheorem{definition}[theorem]
{\hspace{\parindent }D{\scriptsize EFINITION}}
\newtheorem{problem}[theorem]
{\hspace{\parindent }P{\scriptsize ROBLEM}}
\newtheorem{conjecture}[theorem]
{\hspace{\parindent }C{\scriptsize ONJECTURE}}
\newtheorem{example}[theorem]
{\hspace{\parindent }E{\scriptsize XAMPLE}}
\newtheorem{remark}[theorem]
{\hspace{\parindent }R{\scriptsize EMARK}}
\renewcommand{\thetheorem}{\arabic{section}.\arabic{theorem}}
\renewcommand{\theenumi}{(\roman{enumi})}
\renewcommand{\labelenumi}{\theenumi}

\maketitle

\begin{abstract}
In Euclidean geometry the circle of Apollonious is the locus of points in the plane from which two
collinear adjacent segments are perceived as having the same length. In Hyperbolic geometry, the analog of this
locus is an algebraic curve of degree four which can be bounded or  ``unbounded". We study this locus and give a simple description of this curve
using the half-plane model. In the end, we give the motivation of our investigation and calculate the probability that three collinear adjacent segments
can be seen as of the same positive length under some natural assumptions about the setting of the randomness considered.
\end{abstract}

\noindent
\section{Introduction} 
\begin{figure}[h]
\centering
 \includegraphics[scale=0.2]{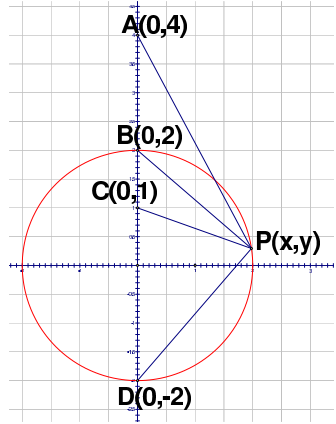}
\caption{The Circle $x^2+y^2=4$ } \label{locus1}
\end{figure}

In most of the textbooks the Circle of Apollonius is discussed in conjunction with the Angle Bisector Theorem:
\emph{``The angle bisector in a triangle divides the opposite sides into a ratio equal to the ratio of the adjacent sides."} Once one realizes that
the statement can be equally applied to the exterior angle bisector, then the  Circle of Apollonius appears naturally (Figure~\ref{locus1}), since the two angle bisectors are perpendicular.

\noindent For instance, an easy exercise in algebra shows that the circle of equation $x^2+y^2=4$ is equivalent to
$$\frac{\sqrt{x^2 +(y-4)^2}}{\sqrt{x^2 +(y-1)^2}}=\frac{BA}{BC}=\frac{DA}{DC}=2,$$

\noindent taking $A(0,4)$, $B(0,2)$, $C(0,1)$ and $D(0,-2)$. Similar calculations can be employed to treat the general situation, i.e., taking
$A(0,a)$, $B(0,b)$ and $C(0,c)$ with real positive numbers $a$, $b$ and $c$ such that $a>b>c$. Then we can state the well-known result:

\begin{theorem}[Apollonius]\label{apollonius}
Given points $A$, $B$ and $C$ as above, the set of points $P(x,y)$ in the plane characterized by the equality $\angle APB\equiv \angle BPC$ is
\par
(i) the line of equation $y=(a+c)/2$ if $b=(a+c)/2$;
\par
(ii) the circle of equation $x^2+y^2=b^2$, if $b<(a+c)/2$ and $b^2=ac$.
\end{theorem}

\noindent Let us observe that the statement of this theorem does not reduce the generality since the coordinate of point $D$ can be shown to be
$y_D=(2ac-bc-ab)/(a+c-2b)$ and one can take the origin of coordinates to be the midpoint of $\overline{BD}$. This turns out to happen precisely  when $b^2=ac$.

What is the equivalent of this result in Hyperbolic Geometry?  We are going to use the half-plane model, $\mathbb H$, to formulate our answer to this question (see Anderson \cite{anderson} for the
terminology and notation used). Without loss of generality we may assume as before that the three points are on the $y$-axis: $A(0,a)$, $B(0,b)$ and $C(0,c)$, with real positive numbers $a$, $b$ and $c$ such that $a>b>c$.

\begin{theorem}\label{theorem1}
Given points $A$, $B$ and $C$ as above, the set of points $P(x,y)$ in the half-plane $\mathbb H$,  characterized by the equality $\angle APB\equiv \angle BPC$ in $\mathbb H$ is
the curve given in polar coordinates by

\begin{equation}\label{eq1}
r^4(2b^2-a^2-c^2)=2r^2(a^2c^2-b^4)\cos (2\theta)+b^2(2a^2c^2-a^2b^2-c^2b^2).
\end{equation}

Moreover, \par
(i) if $b=\left(\frac{a^2+c^2}{2}\right)^{\frac{1}{2}}$, this curve is half of the  hyperbola of equation $$r^2\cos (2\theta)+b^2=0, \ \ \theta\in (\frac{\pi}{4},\frac{3\pi}{4}),$$\par
(ii) if $b=\sqrt{ac}$ this curve is the semi-circle $r=b$, $\theta\in (0,\pi)$,\par
(iii) if $b=\left(\frac{a^{-2}+c^{-2}}{2}\right)^{-\frac{1}{2}}$, this curve is half of the lemniscate
 of equation $$r^2+b^2\cos (2\theta)=0, \ \ \theta\in (\frac{\pi}{4},\frac{3\pi}{4}).$$\par
\end{theorem}

It is interesting that all these particular cases in Theorem~\ref{theorem1}, can be accomplished using integer values of $a$, $b$ and $c$.
This is not surprising for the Diophantine equation $b^2=ac$ since one can play with the prime decomposition of $a$ and $c$ to get $ac$ a perfect square. For the equation
$2b^2=a^2+c^2$ one can take a Pythagorean triple and set $a=|m^2+2mn-n^2|$, $c=|m^2-2mn-n^2|$ and $b=m^2+n^2$ for $m,n\in \mathbb N$. Perhaps it is quite intriguing for
some readers that the last  Diophantine equation $2a^2c^2=a^2b^2+c^2b^2$ is satisfied by the product of some quadratic forms, namely
$$\begin{array}{c} a=(46m^2+24mn+n^2)(74m^2+10mn+n^2),  \\ b=(46m^2+24mn+n^2)(94m^2+4mn-n^2), \\
\text{and}\ \  c=(94m^2+4mn-n^2)(74m^2+10mn+n^2), \ \text{for}\ \ m,n\in \mathbb Z.
\end{array}$$
Another surprising fact is that  Theorem~\ref{theorem1} appears as a more general result than Theorem~\ref{apollonius} because of part (ii) and the observation that we can
always assume that the axes of coordinates are chosen with the origin to be the center of the Apollonius circle.

\begin{figure}[h]
\centering
$\underset{b>\left(\frac{a^2+c^2}{2}\right)^{\frac{1}{2}} }{\includegraphics[scale=0.2]{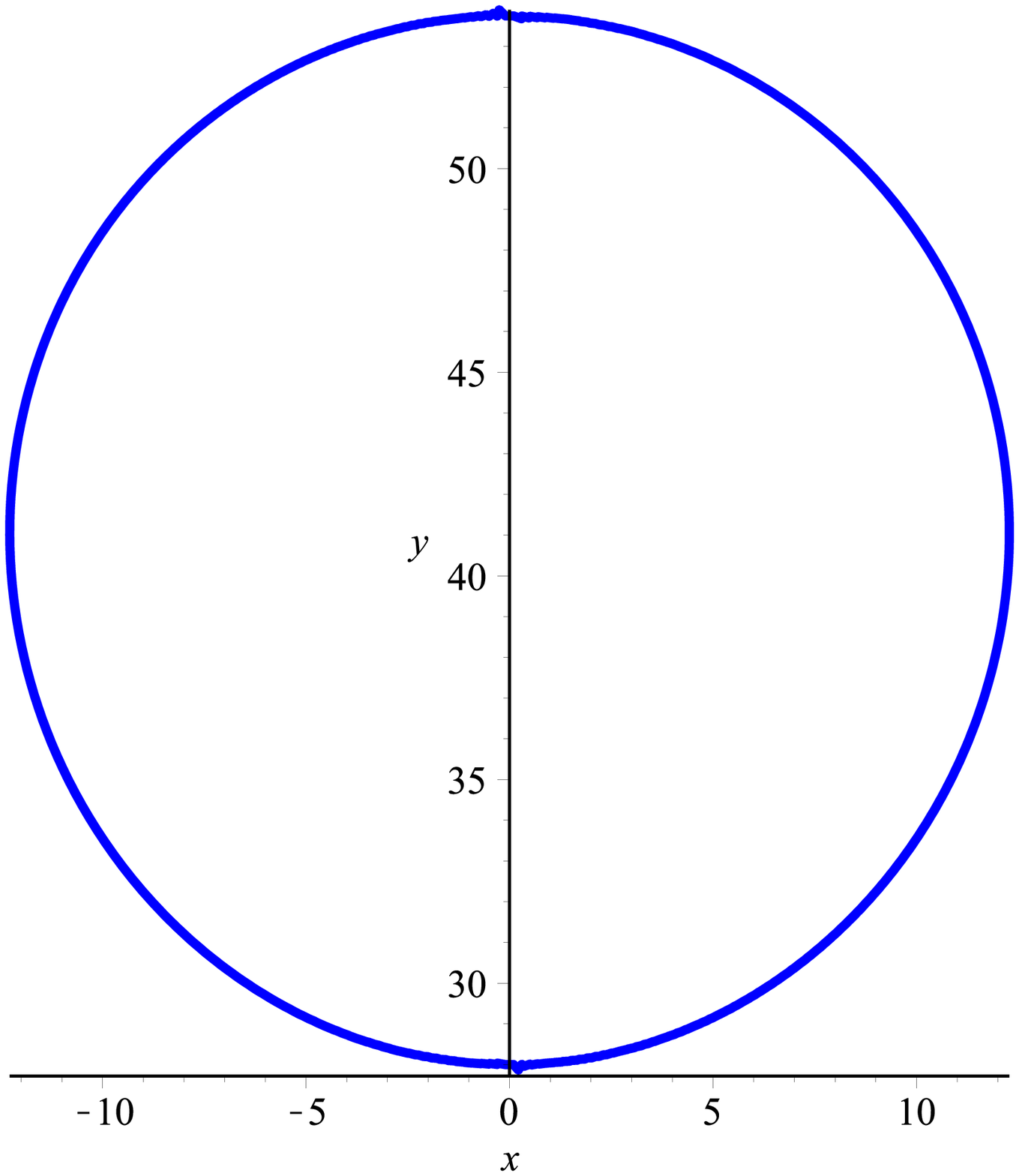}}$
$\underset{b=\left(\frac{a^2+c^2}{2}\right)^{\frac{1}{2}} }{\includegraphics[scale=0.2]{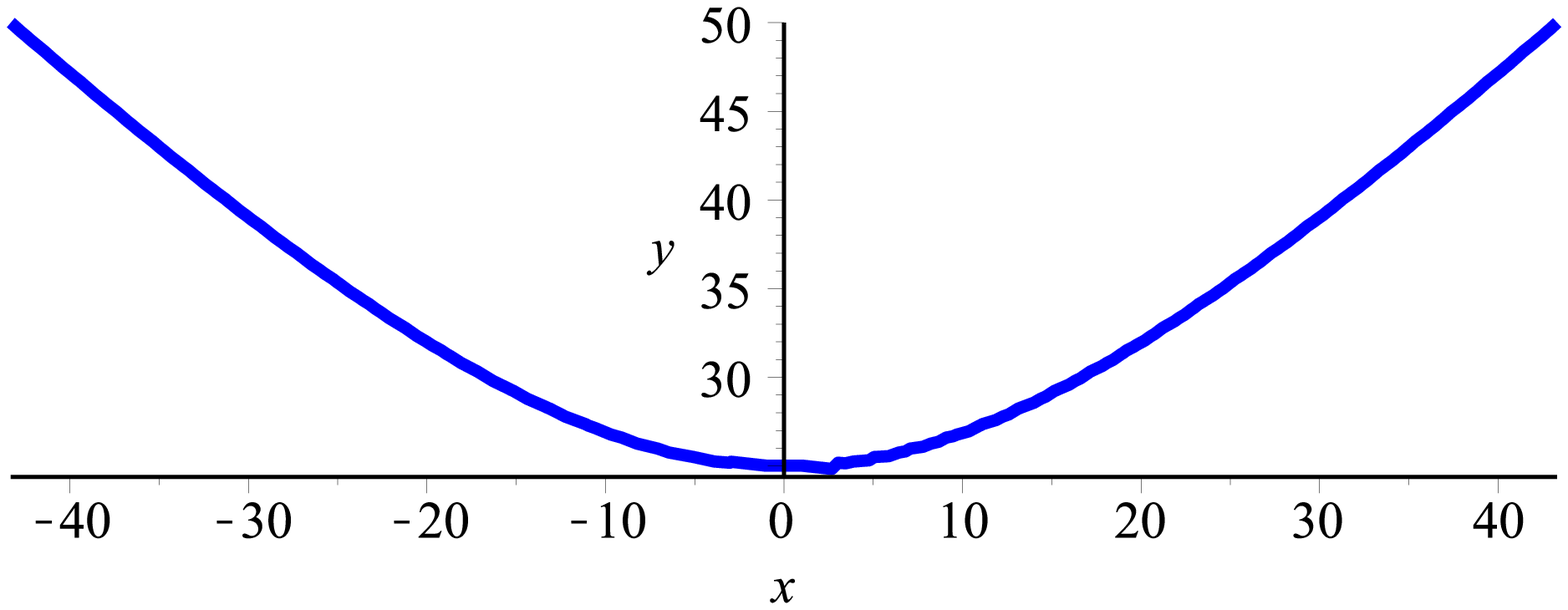}}$
$\underset{\sqrt{ac}<b<\left(\frac{a^2+c^2}{2}\right)^{\frac{1}{2}} }{\includegraphics[scale=0.2]{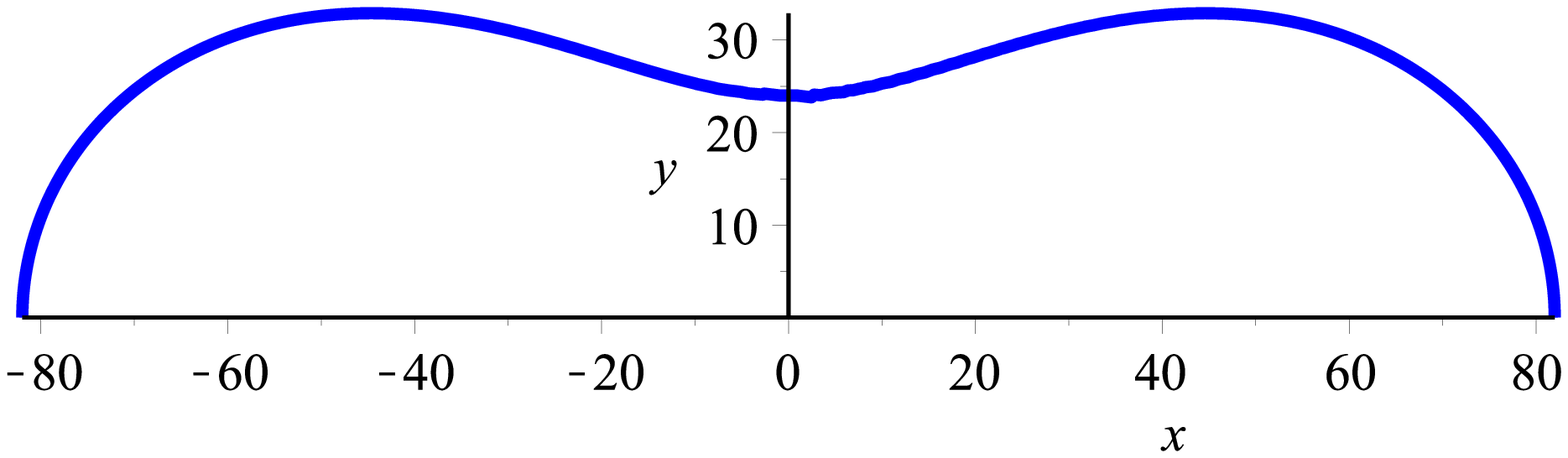}}$
$\underset{b=\sqrt{ac}}{\includegraphics[scale=0.2]{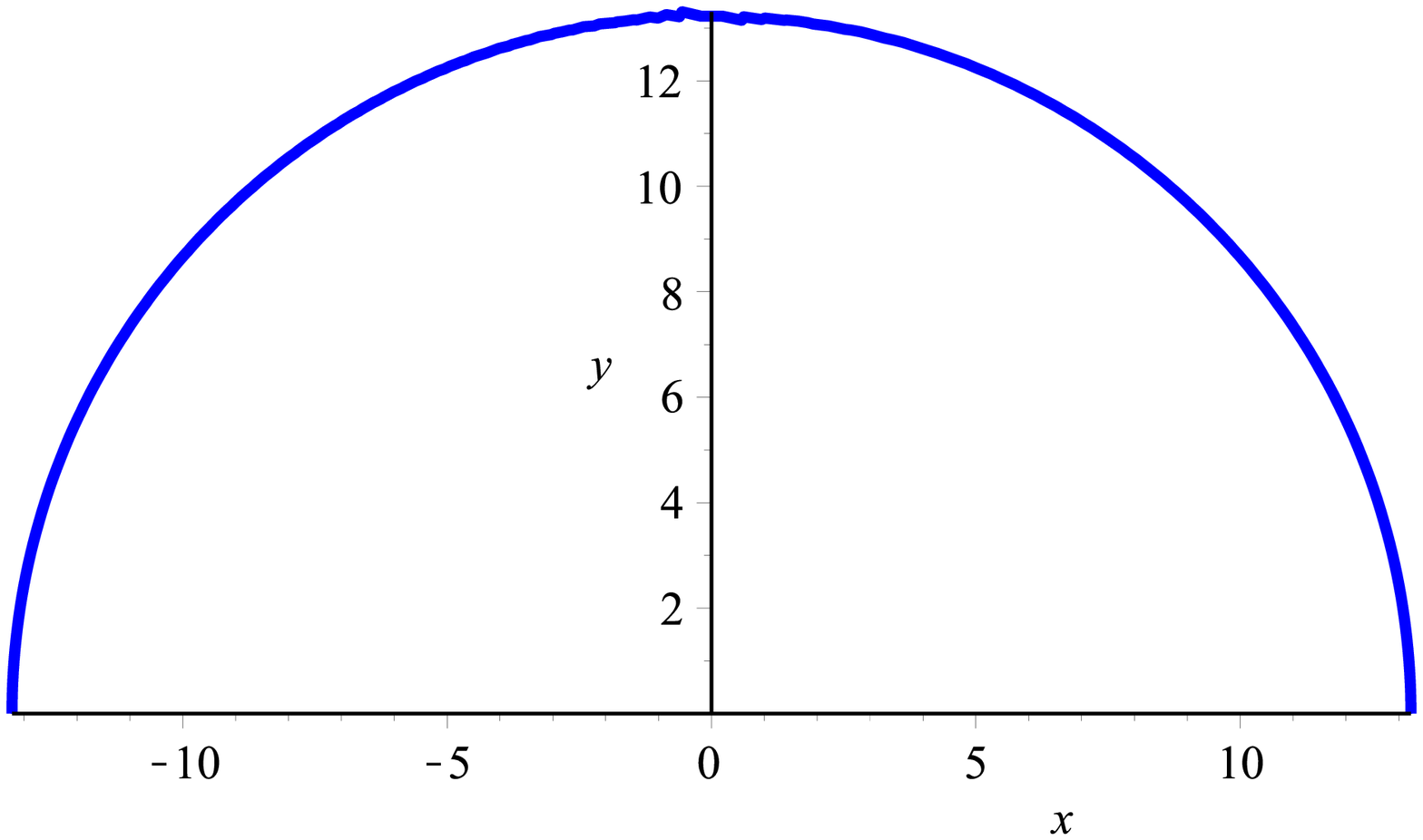}}$
$\underset{\left(\frac{a^{-2}+c^{-2}}{2}\right)^{-\frac{1}{2}}<b<\sqrt{ac}}{\includegraphics[scale=0.2]{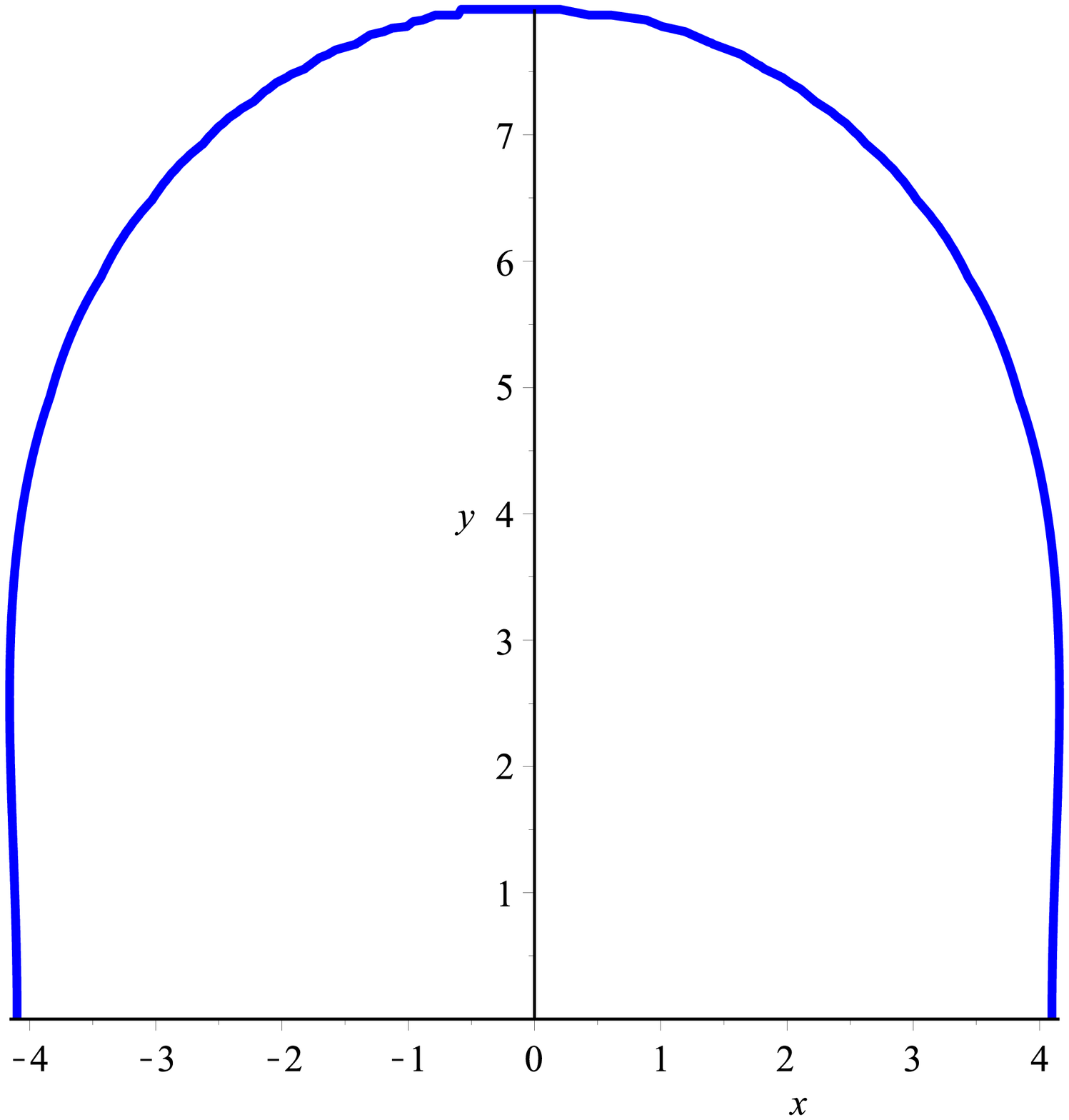}}$
$\underset{b=\left(\frac{a^{-2}+c^{-2}}{2}\right)^{-\frac{1}{2}}}{\includegraphics[scale=0.2]{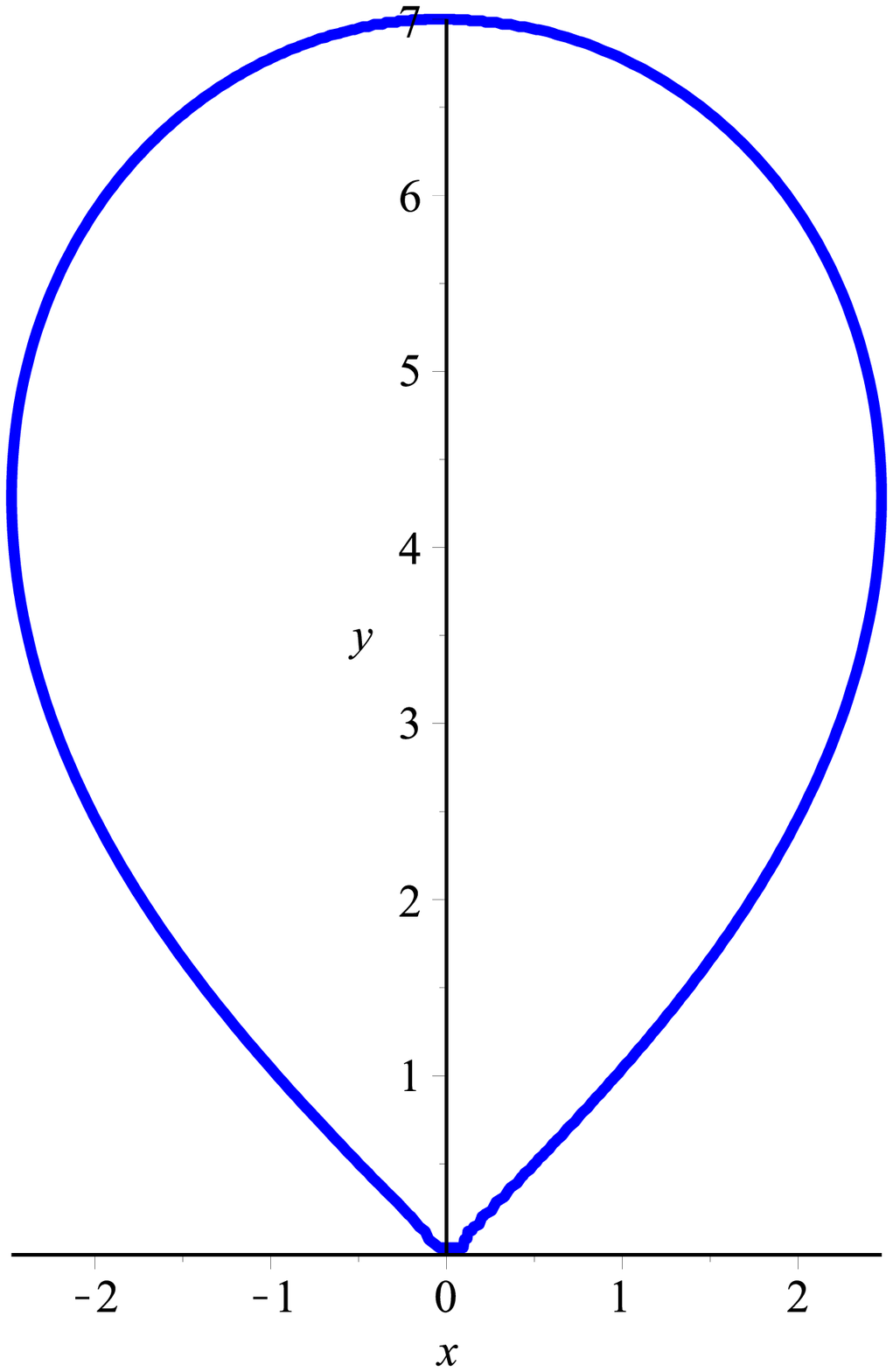}}$
\caption{The curve in polar coordinates (a=35, c=5) } \label{locus2}
\end{figure}

In Figure~\ref{locus2}, we included all of the possible shapes of the locus in Theorem~\ref{theorem1} except for the case
$b<\left(\frac{a^{-2}+c^{-2}}{2}\right)^{-\frac{1}{2}}$ which is similar to the case $b>\left(\frac{a^2+c^2}{2}\right)^{\frac{1}{2}}$. We notice
a certain symmetry of these cases showing that the hyperbola ($b=\left(\frac{a^2+c^2}{2}\right)^{\frac{1}{2}}$) is nothing else
but a lemniscate in hyperbolic geometry. With this identification, it seems like the curves we get, resemble all possible shapes of the intersection of a plane with a torus.

In the next section we will prove the above theorem and in the last section we will give the motivation of our work.
\section{Proof of Theorem~\ref{theorem1}}
\begin{proof} Let us consider a point $P$ of coordinates $(x,y)$  with the given property as in Figure~\ref{locus3} which is not on the line $\overline{AC}$, ($x\not =0$). Then the Hyperbolic lines determined
by $P$ and the three points $A$, $B$ and $C$ are circles orthogonal on the $x$-axis. We denote their centers by $A'(a',0)$, $B'(b',0)$ and $C'(c',0)$.
\begin{figure}[h]
\centering
   \includegraphics[scale=0.3]{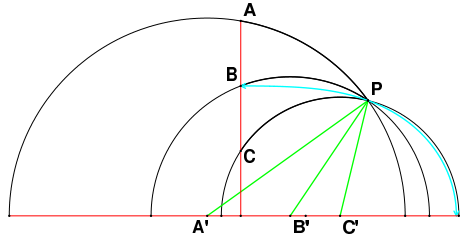}
\caption{The point $P$ and the lines determined by it with $A$, $B$ and $C$} \label{locus3}
\end{figure}
The point $A'$ can be obtained as the intersection of the perpendicular bisector of $\overline{PA}$ and the  $x$-axis.
Similarly we obtain the other two points $B'$ and $C'$. The equation of the perpendicular bisector of $\overline{PA}$ is
$Y-\frac{y+a}{2}=-\frac{x}{y-a}(X-\frac{x}{2})$ and so $a'=\frac{x^2+y^2-a^2}{2x}$. Similar expressions are then obtained for
$b'$ and $c'$, i.e., $b'=\frac{x^2+y^2-b^2}{2x}$ and $c'=\frac{x^2+y^2-c^2}{2x}$. Which shows that the order of the points
 $A'$, $B'$ and $C'$ is reversed ($a'<b'<c'$).  The angle between the Hyperbolic lines $\overset{\leftrightarrow}{PA}$ and $\overset{\leftrightarrow}{PA}$
 is defined by the angle between the tangent lines to the two circles at $P$, which is clearly equal to the angle between the radii corresponding to
 $P$ in each of the two circles. So, $m_{\mathbb H}(\angle APB)=m(\angle A'PB')$ and $m_{\mathbb H}(\angle BPC)=m(\angle B'PC')$.
This equality is characterized by the proportionality given by the Angle Bisector Theorem in the triangle $PA'C'$:
$$\frac{PA'}{PC'}=\frac{A'B'}{B'C'}\Leftrightarrow  \frac{\sqrt{(x^2-y^2+a^2)^2+4x^2y^2}}{\sqrt{(x^2-y^2+c^2)^2+4x^2y^2}}=\frac{a^2-b^2}{b^2-c^2}.$$
Using polar coordinates, $x=r\cos \theta$ and $y=r\sin \theta$, we observe that $x^2-y^2=r^2\cos 2\theta$ and $2xy=r^2\sin 2\theta$.
Hence the above equality is equivalent to $$(r^4+2a^2r^2\cos 2\theta+a^4)(b^2-c^2)^2=(r^4+2c^2r^2\cos 2\theta+c^4)(a^2-b^2)^2.$$
One can check that a factor of $(a^2-c^2)$ can be simplified out and in the end we obtain (\ref{eq1}).\end{proof}

\section{Four points ``equally" spaced and the motivation}
Our interest in this locus was motivated by the Problem 11915 in this Monthly (\cite{kidwell&Meyerson}). This problem stated:
\emph {Given four (distinct) points $A$, $B$, $C$ and $D$
in (this) order on a line in Euclidean space, under what conditions will there be a point $P$ off the line such that the angles
$\angle APB$, $\angle BPC$, and $\angle CPD$ have equal measure?}

It is not difficult to show, using  two Apollonius circles, that the existence  of such a point $P$ is characterized by
the inequality involving the cross-ratio

\begin{equation}\label{eq2}
[A,B;C,D]=\frac{\frac{BC}{BA}}{\frac{DC}{DA}}<3.
\end{equation}
We were interested in finding a similar description for the same question in Hyperbolic space. One can think using the same idea
of the locus that replaces the Apollonius circle in Euclidean geometry and that is why we looked into finding what this locus is.
Having the description of this locus one can see that the intersection of two curves as in Theorem~\ref{theorem1} is difficult to predict.
Fortunately, we can use the calculation done in the proof of Theorem~\ref{theorem1} and formulate a possible answer in the new setting.
Given four points $A$, $B$, $C$ and $D$ in (this) order on a line in the Hyperbolic space, we can use an isometry to transform them on the
line $x=0$ and having coordinates $A(0,a)$, $B(0,b)$, $C(0,c)$ and $D(0,d)$ with $a>b>c>d$. The the existence of a point $P$ off the line $x=0$, such that the angles
$\angle APB$, $\angle BPC$, and $\angle CPD$ have equal measure in the Hyperbolic space is equivalent to the existence of $P$ in Euclidean space corresponding to
the points $A'$, $B'$, $C'$ and $D'$ as constructed in the proof of Theorem~\ref{theorem1}. Therefore the answer is in terms
of a similar inequality

 \begin{equation}\label{eq3}
[A',B';C',D']=\frac{\frac{B'C'}{B'A'}}{\frac{D'C'}{D'A'}}<3 \Leftrightarrow \frac{(b^2-c^2)(a^2-d^2)}{(a^2-b^2)(c^2-d^2)}<3.
\end{equation}

To have a different take of what (\ref{eq2}) means we will translate it into a geometric probability which is not difficult
to compute:  \emph {if two points are randomly selected (uniform distribution) on the segment $\overline{AD}$, then the probability
that a point $P$ off the line $\overline{AD}$ such that the angles
$\angle APB$, $\angle BPC$, and $\angle CPD$ have equal measure exists (where the two points are denoted by $B$ and $C$, $B$ being the closest to $A$),  is equal to $$P_e=\frac{15-16\ln 2}{9}\approx 0.4345$$}

The inequality (\ref{eq3}) gives us the similar probability in the Hyperbolic space:

$$P_h=\frac{2\sqrt{5}\ln(2+\sqrt{5})-5}{5\ln 2}\approx 0.4201514924$$

\noindent where the uniform distribution here means, it is calculated with respect  to the measure $\frac{1}{y}dy$ along the $y$-axis.

\section{Spherical Geometry and the perfect setting}
Our both problems must have an even more interesting answer in the setting of spherical geometry.
Due to the infinite nature of both Euclidean and Hyperbolic spaces the geometric probability question
makes sense only in limiting situations. In this case we can simply ask:
\vspace{0.1in}

{\bf Problem 1:} \emph {What is the equivalent of the circle of Apollonius in spherical geometry?}

\vspace{0.1in}

{\bf Problem 2:} \emph {Given a line in spherical geometry and four points on it, chosen at random with uniform distribution,
what is the probability that the points look equidistant from a point on the sphere that is not on that line?}

\vspace{0.1in}

We find the problems more fascinating because there is no clear order in this geometry.  Given three collinear points, any of them can be thought as between the other two.
So, we may expect three different curves as a result. If three points $A$, $B$ and $C$ are already equidistant on a line $\ell$, the three lines perpendicular to $\ell$ through these points satisfy the 
locus requirement and they are concurrent. Is there such a point (in fact at least two if one exists) of concurrency, for a general position of three points on a line?
Perhaps not, if the three points are close to one another; but there must be quite a variety of situations when it is possible. Also, the probability in question, we anticipate to be somewhat bigger.

\end{document}